\documentclass[12pt]{article}

\usepackage[only,ninrm,elvrm,twlrm,sixrm,sevrm,egtrm,tenrm]{rawfonts}
\usepackage{epsf}
\def\R{I\kern -0,37 em R}
\def\Z{Z\kern -0,60 em Z}
\usepackage{subfigure}
\usepackage{graphicx}
\usepackage{rotating}
\usepackage{morefloats}
\usepackage[T1]{fontenc}
\usepackage[utf8]{inputenc}
\usepackage{amssymb}
\usepackage{amsmath}

\title{The Standard Complex and the 3-dimensional Poincaré Conjecture}

\author{Rui Almeida}

\date{}

\begin{document}

\maketitle

\begin{abstract}
 
We develop a method for constructing standard complexes which turns easy the calculation of their algebraic invariants and, as well, the precise evaluation of whether these complexes are embeddable or not in a 3-manifold. This method applies to all familiar spines of 3-manifolds and, in particular, to the Bing house with two rooms and the classical standard spine of the Poincaré sphere. Finally, we exhibit a compact, connected standard complex which is embeddable into an orientable 3-manifold, its fundamental group is $\Z_{2}$ and it contains a Klein bottle. This standard complex is the spine of a reducible 3-manifold $M^3$, sum of a Seifert fiber space with a fake solid torus, whose universal covering space $W^3$ is a closed and simply connected 3-manifold that cannot be homeomorphic to $S^{3}$.

\end{abstract}


\subsection*{Introduction}

At the end of the second complement to his Analysis Situs, Henri Poincaré claimed that a 3-dimensional manifold with trivial homology is homeomorphic to the standard 3-sphere $S^{3}$. More precisely, he announced in [23], p.308, the following statement: \medskip

\begin{quotation}
{\it ``Tout polyèdre qui a tous ses nombres de BETTI égaux à 1 et tous ses tableaux T{q} bilatères est simplement connexe, c'est-à-dire homéomorphe à l'hypersphère''.}
\end{quotation}

The sentence ``le polyèdre qui a tous ses tableaux T{q} bilatères'' means that all the torsion coefficients of the polyhedron do vanish.\medskip

 Thereafter, in [24], he gave an example of a non-simply connected homology 3-sphere, called nowadays, the Poincaré sphere.\medskip

\begin{quotation}
\it ``On pourrait alors se demander si la considération de ces coefficients suffit; si une variété dont tous les nombres de BETTI et les coefficients de torsion sont égaux à 1 est pour cela simplement connexe au sens propre du mot, c'est-à-dire homéomorphe à l'hypersphère; ou si, au contraire, il est nécessaire, avant d'affirmer qu'une variété est simplement connexe, d'etudier son groupe fondamentale, tel que je l'ai défini dans le $<<$Journal de l'Ecole Polytechnique $>>$ , \S 12, page 60.
 
Nous pouvons maintenant repondre a cette question; 

j'ai formé en effect un exemple d'une variété dont tous les nombres de BETTI et les coefficients de torsion sont égaux à 1, et qui pourtant n'est pas simplement connexe'' ([24], p.46).
\end{quotation} 

We conclude, from the above two assertions that, for Poincaré, the term ``variété simplement connexe'' was understood  as being the standard 3-sphere. Accordingly, as announced on the first page of his {\it Cinquième complément à l'Analyse Situs}, Poincaré carried out the construction of a homology 3-sphere non-homeomorphic to $S^{3}$ and afterwards raised the following question: \medskip

\begin{quotation}
{\it ``Est-il possible que le groupe fondamental de V se réduit à la substituition identique, et que pourtant V ne soit pas simplement connexe?''} ([24], p.110).
\end{quotation}

 It is precisely this question that is known as the Poincaré conjecture.\medskip

It is important, at this point, to mention that the above query was raised immediately after the construction of his counter-example and Poincar\'e certainly did not have, at that time, sufficient data available so as to lead him to conjecture that the only 3-dimensional homotopy sphere was $S^{3}$. Quite to the contrary, he could be led to believe on the existence of further counter-examples.

Finally, Poincaré predicted that the question he had raised would take a very long time to be solved. This prediction merits some further comments.

In his {\it Second complément à l'Analysis Situs}, published in 1900, Poincaré  believed that $S^{3}$ was the only  homology 3-sphere and took almost four years to construct his homology 3-sphere. In fact, he terminated his {\it Cinquième complément à l'Analyse Situs} in november, 1903. Poincaré knew that this question was extremely difficult and that, consequently, an eventual construction of a homotopy 3-sphere non-homeomorphic to $S^{3}$ providing a complete solution to this matter was not to be expected in some near future.\medskip
 
Since then, many efforts have been made in view to solve this ``conjecture'' and one of the first attempts was done by J.H.C. Whitehead in [32]. Finally, we would also like to mention  the following facts:\medskip

a) Last century, many mathematicians have worked in view of solving the Poincaré problem. The natural difficulties, found in dimension 3, put in evidence  many other  related questions and, eventually,  the dimension 3 was extended to an arbitrary dimension. More specifically, the following conjecture was set up ``Let M be an n-dimensional homotopy sphere. Then M is homeomorphic to the standard sphere $S^{n}$.'' The result was proved for all $n \neq 3$ thus reducing the initial conjecture just to the case of $n$ = 3. We remark here that a topological solution to the Poincaré conjecture, in dimension 3, is as well a smooth solution since, in this dimension, the topological and the differentiable categories are equivalent ( see [19] ). This, however, is no longer the case in dimensions larger than 3, where the categories are distinct, and consequently places in evidence the appearance of a smooth Poincaré conjecture with a new type of manifold called an exotic n-sphere namely, an n-manifold homeomorphic but not diffeomorphic to $S^{n}$. It has been proved that there do not exist exotic 5 and 6-spheres, the first examples of such spheres having been exhibited in dimension 7 ([18]). Furthermore, the smooth Poincaré conjecture is still an open problem in dimension 4.
\medskip 

b) Following Bing [1], a closed and connected 3-manifold $M^{3}$ is homeomorphic to $S^3$ if and only if each simple closed curve in $M^{3}$ lies in a 3-cell in $M^3$. We remark, however, that this caracterization of $S^{3}$ is stronger than the simple connectivity.\medskip

 It is today accepted, by the entire mathematical community, that the Poincaré conjecture in dimension 3 is a true statement. Most likely, this opinion will have to change after the appearance of this paper. Our purpose here is precisely to provide an example of a 3-manifold homotopy equivalent to $S^{3}$, called a homotopy 3-sphere, that is not homeomorphic to $S^{3}$. Our approach is based on the theory of simplicial complexes.\medskip

\subsection*{\S1 Basic Theory of Complexes}

We shall assume the basic definitions concerning the theory of complexes as can be found, for example, in the references [9], [12], [27], [29]. Nevertheless, in view of making our text clearer, we shall recall some of the definitions and results that are directly related to this work.\medskip

According to Moise [19] and Bing [2], every topological 3-manifold $M^{3}$ can be triangulated. This means that there exists a pair (${\cal T}$, h), where ${\cal T}$ is a 3-dimensional simplicial complex and  h:$\mid{\cal T}\mid \rightarrow M^{3}$ is a homeomorphism. The underlying set $\mid{\cal T}\mid$ is a polyhedron embedded into $\R^{n}$.\medskip

The above result allows us to define, on $M^{3}$, a piecewise linear structure and, furtermore, piecewise linear maps among the piecewise linear structures and the topological 3-manifolds. We shall abbreviate piecewise linear by p.l.
Thus, if ${\cal K}$ is a 2-dimensional simplicial complex and $M^{3}$ a 3-manifold, we shall say that $\mid{\cal K}\mid$ is p.l.(piecewise linearly) embedded in $M^{3}$, if there is a triangulation (${\cal T}$, h) of $M^{3}$ and a map f:$\mid{\cal K}\mid\rightarrow M^{3}$ such that $h^{-1}$ $\circ$ f is a p.l. embedding. In other terms, $h^{-1}$ $\circ$ f ($\mid{\cal K}\mid$)
  is  a  sub-polyhedron  of $\mid{\cal T}\mid$ ([27], p.7).\medskip

Let ${\cal K}$ be a simplicial complex, ${\cal K}_1$ a sub-complex of ${\cal K}$, $\sigma$ a simplex in ${\cal K}$ and $\tau$ a face of $\sigma$ with dim $\sigma$ = dim $\tau$ + 1. Also, $\tau$ is not a face of any other simplex in ${\cal K}$, meaning that $\tau$ is a free face.\medskip

We say that ${\cal K}_1$ is obtained from ${\cal K}$ by an elementary collapsing if
${\cal K}_1 = {\cal K} -\{\sigma,\tau\}$.

If ${\cal L}$ is a sub-complex of ${\cal K}$ obtained by a finite sequence of elementary collapsings of ${\cal K}$, we shall say that ${\cal K}$ collapses onto ${\cal L}$. It then follows easily, from the definitions, that if ${\cal K}$ collapses onto ${\cal L}$, both complexes have the same Euler characteristic.\medskip

Let P be a polyhedron p.l. embedded in a p.l. n-manifold $M^{n}$. A subset $N$ of $M^n$ is called a regular neighborhood of P in $M^{n}$ if it is a closed neighborhood of P in $M^{n}$, is an n-manifold (with boundary) and  collapses onto P. The following result about existence and uniqueness of regular neighborhoods is sufficient for our purpose.

\medskip

Let P be a compact polyhedron p.l. embedded in a 3-manifold $M^{3}$. Then there exists a regular neighborhood $N$ of P in $M^{3}$ and, moreover, if $N_1$ and $N_2$ are two such regular neighborhoods, there is a p.l.  homeomorphism $f$:$N_{1}\rightarrow N_{2}$ which is the identity on P ([10]).

\medskip

In the literature, the word ``spine'' means, for a closed 3-manifold $M^{3}$, a 2-dimensional sub-complex ${\cal K}$ of $M^{3}$ (for an allowable triangulation of $M^{3}$) such that $M^{3}$ - $\mid{\cal K}\mid$ is a 3-cell (see [17]). If $M^{3}$ is a compact 3-manifold with non-void boundary, ${\cal K}$ is a spine of $M^{3}$ if the latter collapses onto $\mid{\cal K}\mid$ and if, further, there does not exist any elementary collapsing of ${\cal K}$ (see [5],[34]). For a given closed 3-manifold $M^{3}$, a sub-complex ${\cal K}$ and a regular neighborhood $N$ of  $\mid {\cal K} \mid$, we shall say that ${\cal K}$ is a {\it spine} of $M^{3}$ if $M^{3} - Int(N)$ is a disjoint union of 3-balls, where  $Int(N)$ is the interior of $N$. Finally, if $N$ is a compact 3-manifold with boundary, the following relation between the Euler characteristic numbers of $N$ and of $\partial N$ holds:

$$ 2.\chi (N) = \chi (\partial N)$$

\subsection*{\S2 Embedding results in the theory of 3-manifolds}

We now mention two important results that will be useful in proving that our construction provides an authentic counter-example to the 3-dimensional Poincaré conjecture.

 The concept of standard complex, to be developed in the next section, as well as that of fake surface are very closely related. Thus, every standard complex is a closed fake surface and its importance lies in the fact that it is meaningful in deciding whether or not a given complex is orientably 3-thickenable, that is to say, decide whether there exists an embedding into an orientable 3-manifold. In this respect, we quote the following two results that are well adapted to our purpose: \medskip

\begin{quotation}

{\it A fake surface is orientably 3-thickenable if and only if it does not contain a union of the $M\ddot{o}bius$ band with an annulus (one of the boundary circles of the annulus attached to the middle circle of the $M\ddot{o}bius$ band with a map of degree one).}[25] 

\end{quotation}

\begin{quotation}

{\it It is easy to embed a closed regular neighborhood U of the intrinsic 1-skeleton G of K (i.e.: the set of non-manifold points of K - in a (possibly nonorientable) handlebody $\bar H$ so that $Fr(U) \subset \partial \bar H$ and so that G is a spine of $\bar H$. K is obtained from U by attaching connected surfaces $F_{1},... F_{t}$ to Fr(U) along $\partial F_{1}\cup ... \cup \partial F_{t} = U\cap (\overline{K-U})$. Let $\omega_{1} \in H^{1}({\bar H})$ be the orientation class: $\omega_{1}(C)$ is equal to 1 if C pass through nonorientable 1-handles of $\bar H$  an odd number of times, otherwise it is 0. K can be embedded in some 3-manifold if and only if $\omega_{1}(\partial F_{1})$ = ... = $\omega_{1}(\partial F_{t})$ = 0.} [11]

\end{quotation}

\subsection*{\S3 Standard Complexes }

The concept of standard complex is fundamental in order to construct 3-manifolds.
Casler, in [5], proves that each compact 3-manifold with non-void boundary has a standard complex for its spine.

Given a standard complex ${\cal K}$, it is crucial, in general,  to know whether its polyhedron $\mid{\cal K}\mid$ can be embedded into a 3-manifold. When $\mid{\cal K}\mid$ is embeddable then we can take a regular neighborhood $N$ of $\mid {\cal K} \mid$ in this 3-manifold. This neighborhood $N$ is a 3-manifold with boundary and ${\cal K}$ is a spine of $N$. In particular, if ${\cal K}$ is compact, connected and if further its fundamental group is finite and different from $\Z_{2}$, then $N$ is a compact, connected and orientable 3-manifold with finite fundamental group. Its boundary $\partial N$ is a finite union of 2-dimensional spheres and the number of these 2-dimensional spheres is precisely equal to the Euler characteristic of ${\cal K}$. When this fundamental group is equal to $\Z_{2}$, either $N$ is a compact, connected and orientable 3-manifold and $\partial N$ is a finite union of 2-dimensional spheres or else, $N$ is homotopy equivalent to a punctured $RP^{2} \times I $ (see [7]) and $\partial N$ consists of two spaces $RP^{2}$ and some 2-dimensional spheres. \medskip

{\bf Example 3.1 }: Let $m$ and $l$ be, respectively, a meridian and a parallel in a 2-dimensional torus T. We attach a disc $D_1$ along $m$ and a second disc $D_2$ along $l$. In the complement of ${m\cup l}$, all points in this polyhedron are manifold points and the points belonging to ${m\cup l}$ (the non-manifold points) are of two types. The point ${m\cap l}$ has a neighborhood homeomorphic to an euclidean plane with two half planes attached along a pair of perpendicular lines contained in this plane.
All the other points, belonging to ${m\cup l}$, have a neighborhood homeomorphic to an euclidean plane with a half plane attached along a line in this plane. Let ${\cal U}$ be a regular neighborhood of ${m\cup l}$ in this polyhedron and let us cut ${\cal U}$ transversely to ${m}$ as to ${l}$ at two points distinct from the point ${m\cap l}$. We obtain a set that, due to its importance in this work, deserves a more carefull description. \medskip

We take the following points in $\R^{3}$:\medskip

\begin{tabular}{llllllllllll}

  O & = & (0,0,0) & ,   &     $O_1$  & = & (2,0,0)   & ,    &     $O_2$ & = & (0,2,0) & , \\[1ex]
     
 $O_3$ & = & (-2,0,0) & ,    &    $O_4$ & = & (0,-2,0)  & ,   &   $O_5$ & = & (1,1,0) & , \\[1ex]

 $O_6$ & = & (-1,1,0) & ,    &    $O_7$ & = & (-1,-1,0) & ,   &     $O_8$ & = & (1,-1,0) & , \\[1ex] 

 $O_{11}$ & = & (2,-1,0) & ,   &   $O_{12}$ & = & (2,1,0)  & ,  &  $O_{13}$ & = & (2,0,-1) & , \\[1ex]
  
 $O_{21}$ & = & (-1,2,0) & ,  &  $O_{22}$ & = & (1,2,0)  & ,  &   $O_{23}$ & = & (0,2,1) & ,  \\[1ex]
 
 $O_{31}$ & = & (-2,1,0) & ,  &  $O_{32}$ & = & (-2,-1,0) & , &  $O_{33}$ & = & (-2,0,-1) & , \\[1ex]

 $O_{41}$ & = &  (1,-2,0) & ,  &   $O_{42}$ & = & (-1,-2,0) & , & $O_{43}$ & = & (0,-2,1) & , 
 
\end{tabular}
\vspace{2ex}

 and use these points to define, in $\R^{3}$, the following three polygonal lines: \medskip

\begin{tabular}{lp{10.0cm}}
$\alpha_{1}$ = &\small{\mdseries $\overline{O_{11}O_{12}}\cup \overline{O_{12}O_{5}}\cup \overline{O_{5}O_{22}}\cup \overline{O_{22}O_{21}}\cup \overline{O_{21}O_{6}}\cup \overline{O_{6}O_{31}}\cup \overline{O_{31}O_{32}}\cup \overline{O_{32}O_{7}}\cup \overline{O_{7}O_{42}}\cup \overline{O_{42}O_{41}}\cup \overline{O_{41}O_{8}}\cup \overline{O_{8}O_{11}}$,}\\ [2ex]

$\alpha_{2}$ = &\small{\mdseries $\overline{O_{1}O_{13}}\cup \overline{O_{13}O_{33}}\cup \overline{O_{33}O_{3}}\cup \overline{O_{3}O_{1}}$,}\\[2ex]

$\alpha_{3}$ = & \small{\mdseries $\overline{O_{2}O_{23}}\cup \overline{O_{23}O_{43}}\cup \overline{O_{43}O_{4}}\cup \overline{O_{4}O_{2}}$}
\end{tabular}
\vspace{2ex}

where, $\overline{PQ}$ represent the segment $\{ tP + (1-t)Q  | 0 \leq t \leq 1 \}$ for any two points P and Q in $\R^{3}$. These poligonal lines are the boundary curves of three closed 2-cells denoted by $C_{1}$ , $C_{2}$ and$C_{3}$. We remark that $C_{1}$ is contained in the plane z=0, $C_{2}$ in the plane y=o and $C_{3}$ in the plane x=0. The set $C_{1}\cup C_{2}\cup C_{3}$ is our basic neighborhood; it is p.l. homeomorphic to a neighborhood of ${m\cap l}$ as it was described in our initial example. The point O = (0,0,0) is called a vertex point and the set point $C_{1}\cup C_{2}\cup C_{3}$ is called a vertex neighborhood of O.

\begin{figure}[htbp]
\begin{center}
\leavevmode
\hbox{
\epsfxsize=3.5in
\epsffile{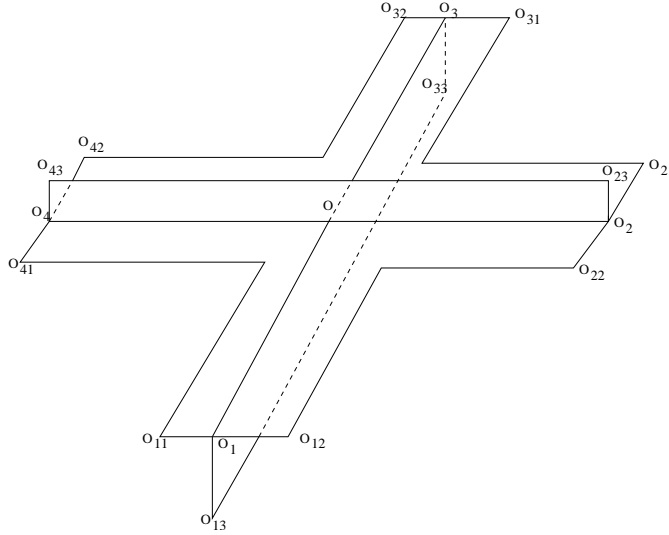}}
\caption{\label{figure 1}This is our vertex neighborhood of O}
\end{center}
\end{figure}

The following subsets of $\partial(C_{1}\cup C_{2}\cup C_{3})$ will be especially important in our oncoming construction of standard complexes: \medskip

\begin{tabular}{llllllll}

$T_{1}$ & = & $\overline{O_{1}O_{11}}\cup \overline{O_{1}O_{12}}\cup \overline{O_{1}O_{13}}$ &, 
$T_{2}$ & = & $\overline{O_{2}O_{21}}\cup \overline{O_{2}O_{22}}\cup \overline{O_{2}O_{23}}$ &,\\[2ex]
$T_{3}$ & = & $\overline{O_{3}O_{31}}\cup \overline{O_{3}O_{32}}\cup \overline{O_{3}O_{33}}$ &,
$T_{4}$ & = & $\overline{O_{4}O_{41}}\cup \overline{O_{4}O_{42}}\cup \overline{O_{4}O_{43}}$ &,

\end{tabular}
\vspace{1ex}

Due to their form, they will be called the four T ends of $C_{1}\cup C_{2}\cup C_3$. \\[0.2ex]

{\bf Definition 3.2 }: Let ${\cal K}$ be a 2-dimensional simplicial complex. Then ${\cal K}$ is a standard complex if, for each non-manifold point P in ${\cal K}$, there exists a neighborhood $V_{P}$ of P in $\mid{\cal K}\mid$ and a p.l. homeomorphism $f : V_{P} \longrightarrow V_{O}$  where f(P) = (0,0,0) and $V_{O}$ is either $C_{1}\cup C_{2}$  or  $C_{1}\cup C_{2}\cup C_{3}$. The sub-complex ${\cal K}^{(1)}$ of ${\cal K}$ that corresponds to the non-manifold points of $\mid{\cal K}\mid$ is called the intrinsic skeleton of ${\cal K}$.\medskip

As a consequence of this definition, let ${\cal K}$ be a compact standard complex, ${\cal K}^{(1)}$ its intrinsic 1-skeleton and ${\cal U}$ a regular neighborhood of $\mid{\cal K}^{(1)}\mid$ in \\ $\mid{\cal K}\mid$. Then, $\mid {\cal K}\mid - {\cal U}$ is a finite and disjoint union of compact surfaces with non-void boundary.
We thus find, in ${\cal U}$, only a finite number of vertex points and consequently ${\cal U}$ is a finite union of vertex neighborhoods of these vertex points glued along their T ends. This suggests that, in order to construct a polyhedron, p.l. isomorphic to $\mid{\cal K}\mid$, we must pick out a finite number of copies of ${C_1}\cup {C_2}\cup {C_3}$, fit all T ends of these neighborhoods appropriately in such a manner as to obtain a polyhedron $\hat {\cal U}$, p.l. isomorphic to ${\cal U}$, and finally attach the above-mentioned finite collection of compact surfaces along the boundary curves of $\hat {\cal U}$.\medskip

To perform this construction, it is very important to have a good description of the matching procedure. 

Let us then take two copies of $C_{1}\cup C_{2}\cup C_{3}$ denoted respectively by $V_A$ and $V_B$, where A and B are the two vertex points.\medskip

We also denote the T ends of $V_A$ by \smallskip

$$T(A)_i = \overline{A_iA_{i1}}\cup \overline{A_iA_{i2}}\cup \overline{A_iA_{i3}}, i = 1, 2, 3, 4$$ \smallskip 

and those of $V_B$ by \smallskip
 
$$T(B)_{j} = \overline{B_{j}B_{j1}}\cup \overline{B_{j}B_{j2}}\cup \overline{B_{j}B_{j3}}, j = 1, 2, 3, 4$$ \smallskip 

Next, we select two indices i and j, a permutation $\sigma : \{1,2,3\} \longrightarrow \{1,2,3\}$ and fit together $\overline{A_iA_{i1}}$ with $\overline{B_{j}B_{j\sigma(1)}}$, $\overline{A_iA_{i2}}$ with $\overline{B_{j}B_{j\sigma(2)}}$, and $\overline{A_iA_{i3}}$ with $\overline{B_{j}B_{j\sigma(3)}}$.

When the two T ends $T(A)_{i}$ and $T(B)_{j}$ are appropriately matched, we denote this matching simply by:
$(A_{i1} A_{i2} A_{i3}) \sim (B_{j\sigma (1)} B_{j\sigma(2)} B_{j\sigma(3)})$.

It is convenient here to exhibit the six matching possibilities given by the permutation group of \{1,2,3\}.

We have the odd matches:

$$(A_{i1} A_{i2} A_{i3}) \sim (B_{j2} B_{j1} B_{j3})$$
$$(A_{i1} A_{i2} A_{i3}) \sim (B_{j3} B_{j2} B_{j1})$$
$$(A_{i1} A_{i2} A_{i3}) \sim (B_{j1} B_{j3} B_{j2})$$\smallskip

and the even matches:

$$(A_{i1} A_{i2} A_{i3}) \sim (B_{j1} B_{j2} B_{j3})$$
$$(A_{i1} A_{i2} A_{i3}) \sim (B_{j3} B_{j1} B_{j2})$$
$$(A_{i1} A_{i2} A_{i3}) \sim (B_{j2} B_{j3} B_{j1})$$\smallskip

As we shall see further, the existence of even matches in the construction of a standard complex can produce obstructions to the embeddability of this complex into a 3-manifold.\medskip

We next present a practical condition that permit us to know exactly whether a given standard complex embeds or not in some orientable 3-manifold. Let ${\cal K}$ be a compact standard complex, $\mid{\cal K}^{(1)}\mid$ its intrinsic 1-skeleton and ${\cal U}$ a regular neighborhood of $\mid{\cal K}^{(1)}\mid$ in $\mid{\cal K}\mid$. As we have seen previously, ${\cal U}$ = ${\bigcup\limits_{i\in I}(V_{A_i}/ \sim})$, where each $A_{i}$ is a vertex point, $V_{A_{i}}$ is the vertex neighborhood of $A_{i}$ and ``$\sim$'' means that all T ends of these vertex neighborhoods are appropriately glued. Let us take two of these vertex points $A_{i}$ and $A_{j}$ connected by an open 1-cell $a_{l} \subset \mid{\cal K}^{(1)}\mid$. Then, one T end of $V_{A_{i}}$ is matched with one T end of $V_{A_{j}}$  to form $a_{l}$. We say that $a_{l}$ is {\itshape {even}} if the above-mentioned matching is even. Now, each oriented closed path $\gamma \subset Fr({\cal U})$ is homotopic, in ${\cal U}$, to a closed path $\bar a_{1}.\bar a_{2} ... \bar a_{n}$, where each $a_{i}$ is an open oriented 1-cell contained in $\mid{\cal K}^{(1)}\mid$ having exactly two vertex points at its ends and $\bar a_{i}$ means $a_{i}$ together with its two vertex points. We observe that if any of the 1-cells forming the closed path $\bar a_{1}.\bar a_{2}... \bar a_{n}$ is even, then the part of the handlebody where ${\cal U}$ is embedded and that contains $\bar a_{1}.\bar a_{2}... \bar a_{n}$ is a non-orientable handle.\medskip

Having settled these preliminaires, we can rewrite the  Kranjc remark 1 ( [11], p.310 ) in the following manner:\medskip

{\sl Let ${\cal K}$ be a standard complex, ${\cal K}^{(1)}$ its intrinsic skeleton and ${\cal U}$ a regular neighborhood of $\mid{\cal K}^{(1)}\mid$ in $\mid{\cal K}\mid$. 
Then $\mid{\cal K}\mid$ can be embedded into an orientable 3-manifold if and only if, for each connected component $\gamma$ of $Fr({\cal U})$, the closed path $\bar a_{1}.\bar a_{2}. ... \bar a_{n}$, which is homotopic to $\gamma$ in ${\cal U}$, has an even number of even 1-cells.
In particular, the non existence of even matchings in ${\cal U}$ implies that ${\cal U}$ is embeddable in an orientable handlebody and consequently there is no obstruction to embed $\mid{\cal K}\mid$ in some 3-manifold.}\medskip

 {\bf Remark:} If $\mid{\cal K}\mid$ satisfies the above-mentioned condition, then  $\mid{\cal K}\mid$ does not contain a $M\ddot{o}bius$ band with an annulus, where one of the boundary  circles of the annulus is attached to the middle circle of the $M\ddot{o}bius$ band. By Corollary 1.2 of [25], $\mid{\cal K}\mid$ is orientably 3-thickenable. This reinforces our criterion.

 We do not pursue the above discussion in the case where there exist even matchings in ${\cal U}$, but  provide a simple and very suggestive example.
\\[0.2ex]

{\bf Example 3.3 }: Let $V_0 = C\cup \hat C$, where \smallskip

$$C = \{ (x, y, 0) | (x, y) \in [-2, 2]\times [-1, 1]\}$$ 
and 
$$\hat C = \{(x, 0, z) | x \in [ -2, 2] , z \in [0, 1] \}$$.

\begin{figure}[htbp]
\begin{center}
\leavevmode
\hbox{%
\epsfxsize=4.0in
\epsffile{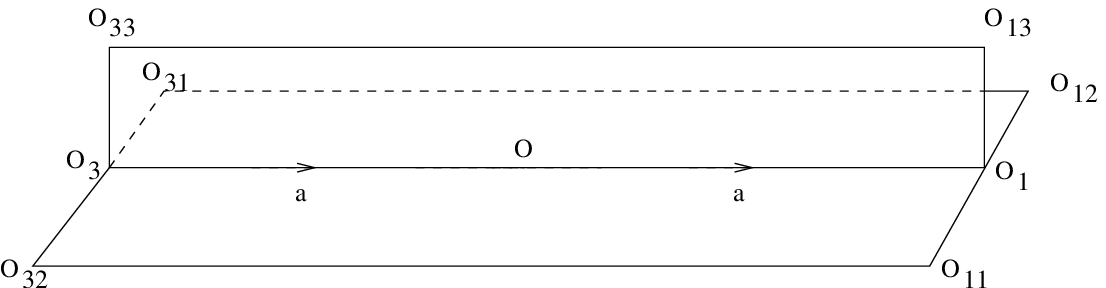}}
\caption{\label{}}
\end{center}
\end{figure}

$V_{0}$ has only the two following T ends: \medskip

$$T(V_{0})_{1} = \overline{O_{1}O_{11}}\cup \overline{O_{1}O_{12}}\cup \overline{O_{1}O_{13}}$$ 

and

$$T(V_{0})_{3} = \overline{O_{3}O_{31}}\cup \overline{O_{3}O_{32}}\cup \overline{O_{3}O_{33}}$$

where

$$O_{1} = (2, 0, 0) , O_{3} = (-2, 0, 0) , O_{11} = (2, -1, 0) , O_{12} = (2, 1, 0) ,$$

$$O_{13} = (2, 0, 1) , O_{31} = (-2, 1, 0) , O_{32} = (-2, -1, 0) , O_{33} = (-2, 0, 1).$$

We match these two T ends by the even relation

$$(O_{11} O_{12} O_{13}) \sim (O_{31} O_{32} O_{33})$$.\medskip

It is easy to see that the non-manifold points of the quotient space $V_{0}/\sim$ is a circle a and $Fr(V_{0}/\sim)$ is the disjoint union of two closed paths $\gamma_{1}$ and $\gamma_{2}$. The path $\gamma_{1}$ is homotopic to a in $V_{0}/\sim$ and $\gamma_{2}$ is homotopic to  a.a in $V_{0}/\sim$. According to our embedding conditions, the polyhedron $\mid{\cal K}\mid$ obtained from $V_{0}/\sim$ by attaching two disks $D_{1}$ and $D_{2}$ along $\gamma_{1}$ and $\gamma_{2}$ respectively, is not embeddable into a 3-manifold.\smallskip

For a more direct description of $\mid{\cal K}\mid$, we can take the projective real plane $RP^{2}$ and attach a disk along a line $a \subset RP^{2}$.\smallskip

Obviously, $\mid{\cal K}\mid$ is compact, connected and simply connected.\smallskip

{\bf Remark}: If we match the two ends by the odd matching relation 
$$(O_{11} O_{12} O_{13}) \sim (O_{31} O_{33} O_{32}),$$
we find that $Fr(V_{O}/\sim)$ is a connected closed path. Attaching a disk along this path, we obtain the classical standard spine of the lenticular space $L_{3,1}$ (or $L_{3,2}$). \medskip

\subsection*{\S4 The construction procedure}

As pointed out in the preceding paragraph, a standard complex has a kind of a decomposition where the main piece is a vertex neighborhood. Using such neighborhoods together with the previously defined procedure of matching together their T ends, we can construct many different examples of standard complexes. Our construction renders easy the determination of their fundamental groups and consequently the determination of the most important invariants of these complexes. This construction obeys the following steps:\smallskip

a) We select n copies $V_{A_i}$, i = 1, ..., n of $C_{1}\cup C_{2}\cup C_{3}$, labeled by their vertex points.\smallskip

b) The T ends, 4n in number, are separated in 2n pairs and, for each such pair, we determine a matching relation.\smallskip

c) When all the T ends are appropriately matched, we set  ${\cal U}$ = ${\bigcup\limits_{i = 1}^{i = n}(V_{A_i}/ \sim})$ and observe  that $Fr({\cal U})$ is composed by a finite number of closed paths $\gamma_{1}, ..., \gamma_{t}$ where we can subsequently attach t disks $D_{1}$, ..., $D_{t}$ to ${\cal U}$ along $\gamma_{1}$, ..., $\gamma_{t}$ respectively, thus obtaining a polyhedron $\mid{\cal K}\mid$, where ${\cal K}$ is a standard complex.\smallskip

It is easy to see that ${\cal K}$ is compact and connected and its Euler characteristic is $n - 2n + t$. The next step consists in determining the fundamental group of ${\cal K}$.\smallskip

d) We next label each of the 2n 1-cells that, jointly with the vertex points, constitutes $\mid{\cal K}^{(1)}\mid$. The attachment, thereafter, of t disks $D_{1}$, ..., $D_{t}$ creates t homotopy relations.\smallskip

Finally, we select one of the vertex points as the base point and take auxiliary paths provided by $(n - 1)$ paths in $\mid{\cal K}\mid$ connecting each of the remaining vertex points to the base point. Defining the representatives of the generators of the fundamental group of $\mid{\cal K}\mid$ and using the homotopy relation obtained in d), we can determine this group. The first homology group of $\mid{\cal K}\mid$ is simply its abelianized group. It is also easy to calculate the Betti numbers $b_{0}$ , $b_{1}$ and $b_{2}$ of $\mid{\cal K}\mid$.\smallskip

{\bf Remark}: The standard complexes described above are called special complexes. More generally, it is possible, in c), to attach some surfaces with non-void boundaries in the place of disks. In this case, we shall change some homotopy relations by homology relations. \medskip

\subsection*{\S5 Some Standard Complexes}
We now give some examples of standard complexes. The Bing house with two rooms and the classical standard spine of the Poincar\'e sphere are two noteworthy examples. The simplest are obtained by just selecting a single vertex neighborhood for which we have 108 possibilities for matching their four T ends but, due the existence of symmetries, we just obtain 14 non-homeomorphic standard complexes and just 4 of these are embeddable into a 3-manifold. As we shall see in the sequel, one of these simplexes is simply connected, its Euler characteristic is 1 and it is embeddable in $\R^{3}$. Moreover, it is a standard spine of a 3-ball. Comparing this complex with the Bing house with two rooms (the classical standard spine of a 3-ball) we realize that this new complex is simpler. More generally, when we select two or more vertex neighborhoods, the task in describing and classifying all the possibilities become considerably more arduous. \medskip

{\bf Remark}: In the next figures, the points $O_{11}$, ...,$O_{43}$, $A_{11}$, ..., $A_{43}$, as well all the other points that define the T ends of the vertex neighborhoods are missing. For any doubt concerning their actual position on the vertex neigborhood, confer the figure 1. \\[0.2ex]

{\bf Example 5.1 }:
 We select just one vertex neighborhood $V_{O}$ as is shown in the figure below and proceed as follows:

\begin{figure}[htbp]
\begin{center}
\leavevmode
\hbox{%
\epsfxsize=3.5in
\epsffile{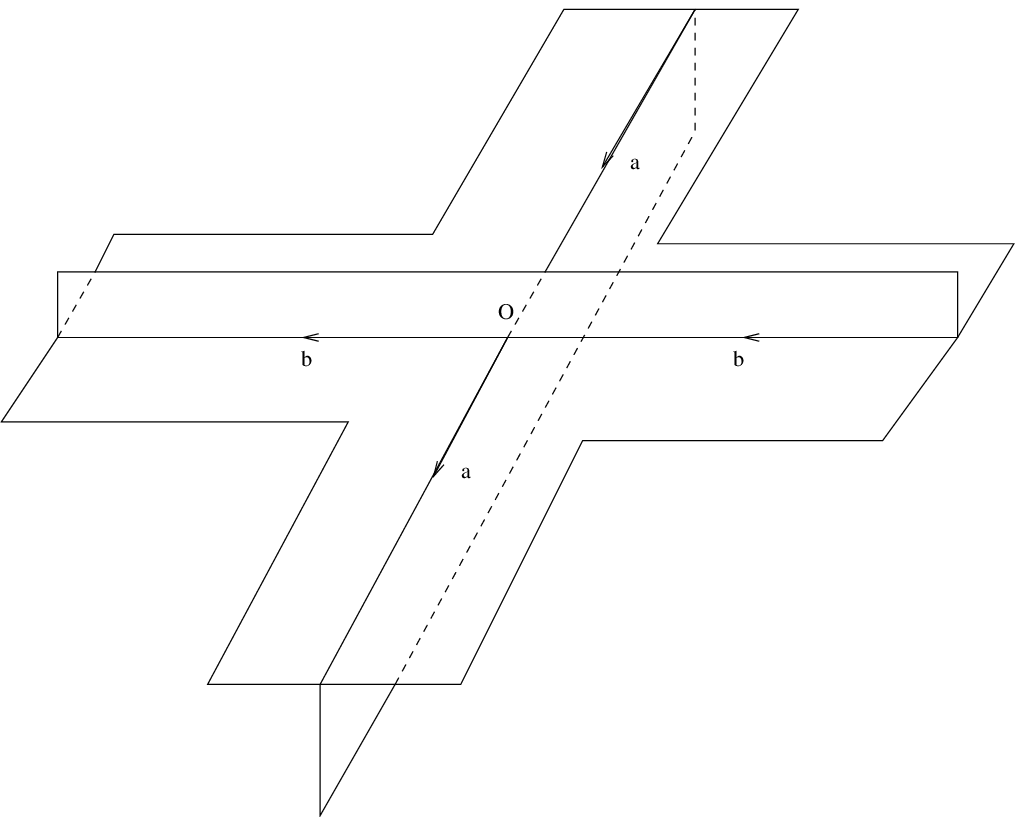}}
\caption{\label{}}
\end{center}
\end{figure}

a) We match the T ends of $V_{O}$ according to the relations:

$$(O_{11} O_{12} O_{13}) \sim (O_{32} O_{31} O_{33})$$

and 

$$(O_{21} O_{22} O_{23} \sim (O_{41} O_{43} O_{42})$$. \medskip

 $Fr(V_{O}/\sim)$ consists of two closed paths. We attach two disks along these paths, getting $\mid{\cal K}\mid$ and the following homotopy relations: \medskip

$$1) \hspace{0,5cm}  ab^{-1}ab^{2} \sim 1 \hspace{4cm}     2) \hspace{0,5cm} a \sim 1$$ \medskip

Let O be the base point. Then A = a and B = b are representatives of the generators of the fundamental group $\pi_{1}(\mid{\cal K}\mid, O)$ satisfying \medskip

$$1) \hspace{0,5cm} AB^{-1}AB^{2} \sim 1 \hspace{4cm}     2) \hspace{0,5cm} A \sim 1$$ \medskip

Thus, $\mid{\cal K}\mid$ is compact, connected and simply connected and it embeds in $\R^{3}$,  a regular neighborhood of $\mid{\cal K}\mid$ in $\R^{3}$ being a 3-ball.\medskip

b) The matching relation: 
$$(O_{11} O_{12} O_{13}) \sim (O_{32} O_{31} O_{33})$$

and

$$(O_{21} O_{22} O_{23}) \sim (O_{42} O_{41} O_{43})$$

produce a 2-dimensional torus $T^{2}$ with two disks attached along a meridian and a parallel of this torus.\medskip

Our standard complex is simply connected and has Euler characteristic equal to 2. A regular neighborhood $N$ has its boundary composed by spheres of dimension 2 and, in gluing two 3-balls to $N$, we obtain $S^{3}$. \medskip

These are the standard complexes where the matching relations are odd and thus produce complexes embeddable into orientable 3-manifolds.\smallskip
The following complex is defined by even matchings:\medskip

c) The matching relations:
$$(O_{11} O_{12} O_{13}) \sim (O_{31} O_{32} O_{33})$$ 

 and 
 
$$(O_{21} O_{22} O_{23})  \sim (O_{41} O_{42} O_{43})$$ 

produce a complex $\mid{\cal K}\mid$ isomorphic to a real projective plane with two disks attached along two different lines of this plane. Obviously, $\mid{\cal K}\mid$ is compact, connected and simply connected, but it cannot be embedded into a 3-manifold.\\[0.2ex]

{\bf Example 5.2}: We select two vertex neighborhoods $V_{A}$ and $V_{B}$

\begin{figure}[htbp]
\begin{center}
\leavevmode
\hbox{%
\epsfxsize=6.0in
\epsffile{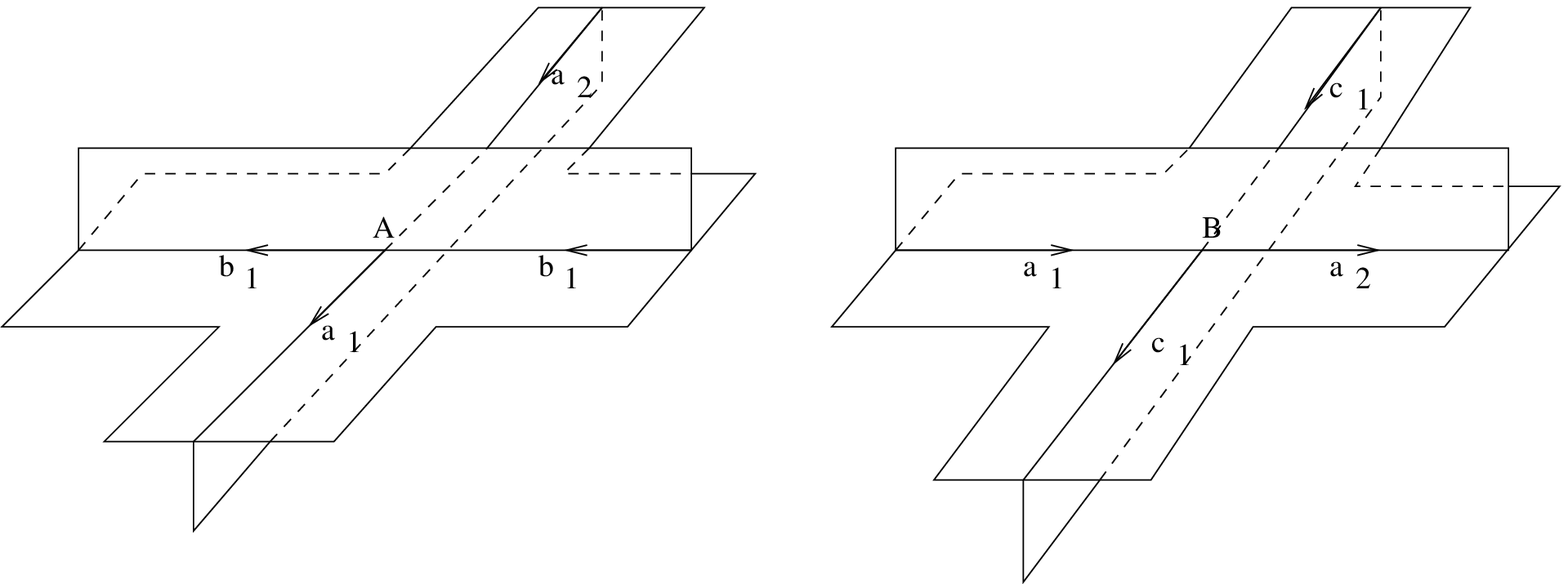}}
\caption{\label{}}
\end{center}
\end{figure}

 and match the T ends of $V_{A}$ and $V_{B}$ by the following relations:\medskip

\begin{tabular}{lllllll}
$(A_{11} A_{12} A_{13})$ & $\sim$ & $(B_{43} B_{42} B_{41})$ & , & $(A_{21} 
A_{22} A_{23})$ & $\sim$ & $(A_{42} A_{41} A_{43})$  \\[1.ex]

$(A_{31} A_{32} A_{33})$ & $\sim$ & $(B_{21} B_{23} B_{22})$ & , & $(B_{11} B_{12} B_{13})$ & $\sim$ & $(B_{32} B_{31} A_{33})$
\end{tabular}
\\[0.2ex]

Then, $Fr(V_{A}\cup V_{B} /\sim)$ consists of three closed paths and, attaching disks along them, we obtain a complex $\mid{\cal K}\mid$ that satisfies the following homotopy relations:\medskip

1) $a_{1} c_{1} a_{1}^{-1} b_{1}^{-1} a_{1} a_{2} b_{1} a_{2}^{-1} c_{1}^{-1}  a_{2} \sim 1$\smallskip

2) $b_{1} \sim 1$\smallskip

3) $c_{1} \sim 1$ \medskip

Let A be the base point and $a_{1}$ the auxiliary path. Then, $A_{1} = a_{1} a_{1}^{-1}$ , $A_{2} = a_{1} a_{2}$ ,  $B_{1} = b_{1}$  ,  $C_{1} = a_{1} c_{1} a_{1}^{-1}$  are representatives of the generators of $\pi_{1}(\mid{\cal K}\mid, A)$ satisfying:\smallskip

1) $C_{1} B_{1}^{-1} A_{2} B_{1} A_{2}^{-1} C_{1}^{-1} A_{2}$ \smallskip
    
2) $B_{1} \sim 1$ \smallskip 
  
3) $C_{1} \sim 1$.\medskip

It is easy to show that $\mid{\cal K}\mid$ is simply connected, its Euler characteristic is 1 and $\mid{\cal K}\mid$ is a standard spine of a 3-ball.\smallskip

In fact $\mid{\cal K}\mid$ is the Bing house with two rooms.
\\[0.2ex]

{\bf Example 5.3}: We select five vertex neighborhoods as shown in the next figure.

\begin{figure}[htbp]
\begin{center}
\leavevmode
\hbox{%
\epsfxsize=6.5in
\epsffile{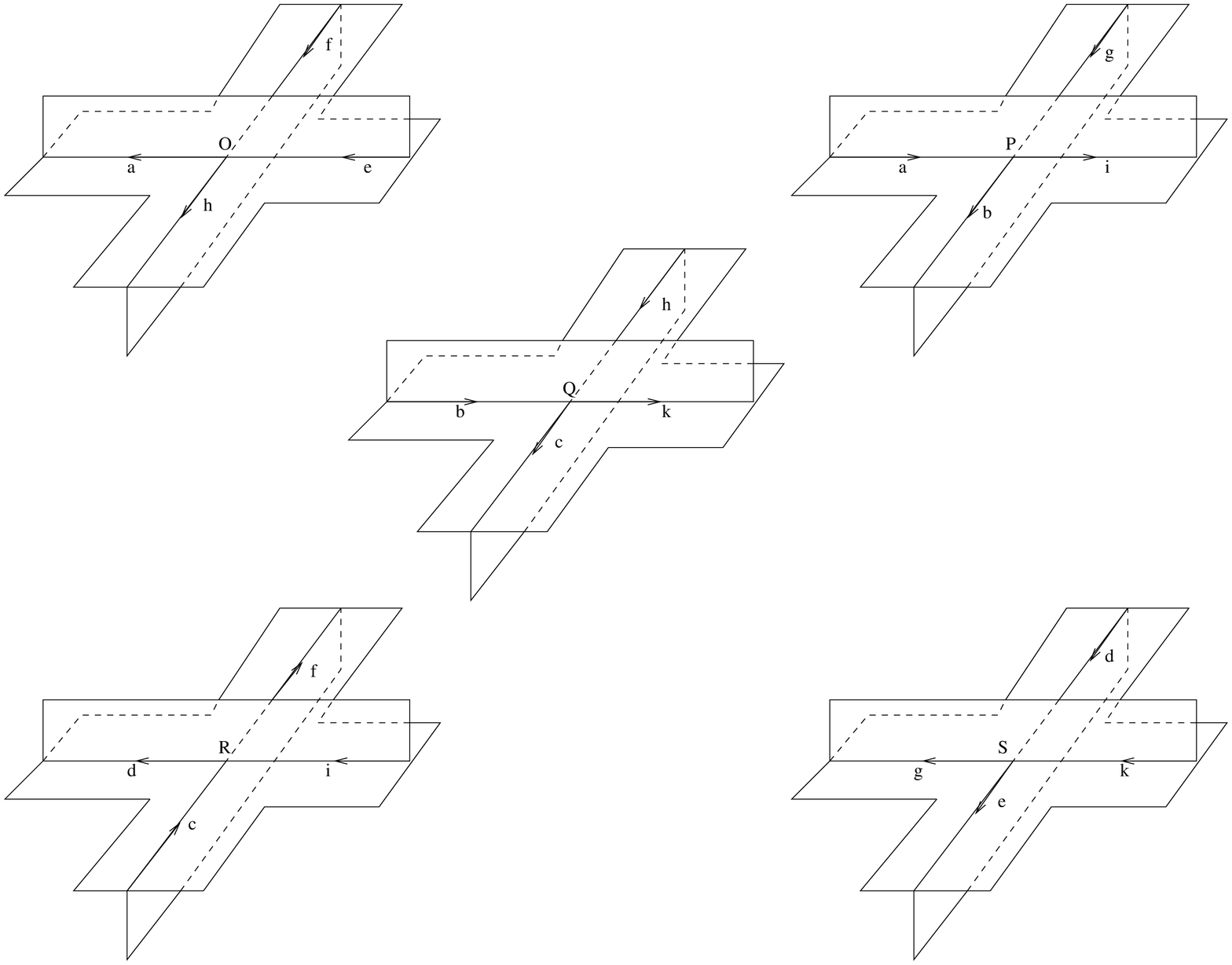}}
\caption{\label{}}
\end{center}
\end{figure}

The matching relations are given by:\medskip

\begin{tabular}{lllllll}
$(O_{11} O_{12} O_{13})$ & $\sim$ & $(Q_{31} Q_{33} Q_{32})$ & , & $(O_{21} O_{22} O_{23})$ & $\sim$ & $(S_{12} S_{11} S_{13})$  \\[1.0ex]

$(O_{31} O_{32} O_{33})$ & $\sim$ & $(R_{31} R_{33} R_{32})$ & , & $(O_{41} O_{42} O_{43})$ & $\sim$ & $(P_{43} P_{42} P_{41})$  \\[1.0ex]

$(P_{11} P_{12} P_{13})$ & $\sim$ & $(Q_{41} Q_{43} Q_{42})$ & , & $(P_{21} P_{22} P_{23})$ & $\sim$ & $(R_{22} R_{21} R_{23})$  \\[1.0ex]

$(P_{31} P_{32} P_{33})$ & $\sim$ & $(S_{41} S_{43} S_{42})$ & , & $(Q_{11} Q_{12} Q_{13})$ & $\sim$ & $(R_{11} R_{13} R_{12})$  \\[1.0ex]

$(Q_{21} Q_{22} Q_{23})$ & $\sim$ & $(S_{21} S_{23} S_{22})$ & , & $(R_{41} R_{42} R_{43})$ & $\sim$ & $(S_{33} S_{32} S_{31})$  
\end{tabular}
\\[0.1ex]

Let ${\cal U}$ be the quotient space $V_{O} \cup V_{P} \cup V_{Q} \cup V_{R} \cup V_{S} /\sim$. Then $Fr({\cal U})$ consists of six closed paths. We attach disks along these closed paths and obtain a polyhedron $\mid{\cal K}\mid$ with the following homotopy relations: \smallskip

\begin{tabular}{lllll}
1) $a b c d e \sim 1$     & &    2) $b k e f^{-1} i^{-1} \sim 1$     & &     3) $a i d k^{-1} h^{-1} \sim 1$ \\[1ex]

4) $c i^{-1} g^{-1} e h \sim 1$     & &     5) $b h^{-1} f^{-1} d g \sim 1$     & &     6) $a g^{-1} k^{-1} c^{-1} f \sim 1$ \medskip 
\end{tabular}

Taking O as the base point and a, h, $f^{-1}$, $f^{-1}d$ as auxiliary paths leading to the vertices P, Q, R and S respectively, the generating path classes of the fundamental group of $\mid{\cal K}\mid$ are represented by: \medskip

$A = a a^{-1}$ , $B = a b h^{-1}$ , $C = h c f$ , $D = f^{-1}d(d^{-1}f)$ , $E = (f^{-1}d)e$ , \smallskip
 
$F = f^{-1} f$ , $G = (f^{-1} d)g a^{-1}$ , $H = h h^{-1}$ , $J = a i f$ , $K = h k(d^{-1}f)$.
\medskip

Then $A \sim D \sim F \sim H \sim 1$ and the homotopy relation can be re-stated in the following form: \medskip

\begin{tabular}{lllll}
1) $B C E \sim 1$  &  &   2) $B K E J^{-1} \sim 1$  &  &   3) $J K^{-1} \sim 1$  \\[1ex] 

4) $C J^{-1} G^{-1} E \sim 1$  &  &   5) $B G \sim 1$  &  &   6) $G^{-1} K^{-1} C \sim 1$ \medskip
\end{tabular} 

An easy calculation, as is done for instance in [29], shows that $\pi_{1}(\mid{\cal K}\mid, O)$ is a group of order 120, whose abelianized group is trivial and the Euler characteristic of $\mid{\cal K}\mid$ is 1. We remark that $\mid{\cal K}\mid$ is the classical standard spine of the Poincar\'e sphere.\\[0.1ex]

The next standard complex has two vertex neighborhoods and its universal covering space has four vertex neighborhoods becoming the spine of a homotopy 3-sphere non-homeomorphic to $S^{3}$.
\\[0.2ex]

{\bf Example 5.4}: We select two vertex neighborhoods $V_{A}$ and $V_{B}$ as shown in the next figure and consider the below indicated matching relations:\\[0.3ex]

\begin{figure}[htbp]
\begin{center}
\leavevmode
\hbox{%
\epsfxsize=6.0in
\epsffile{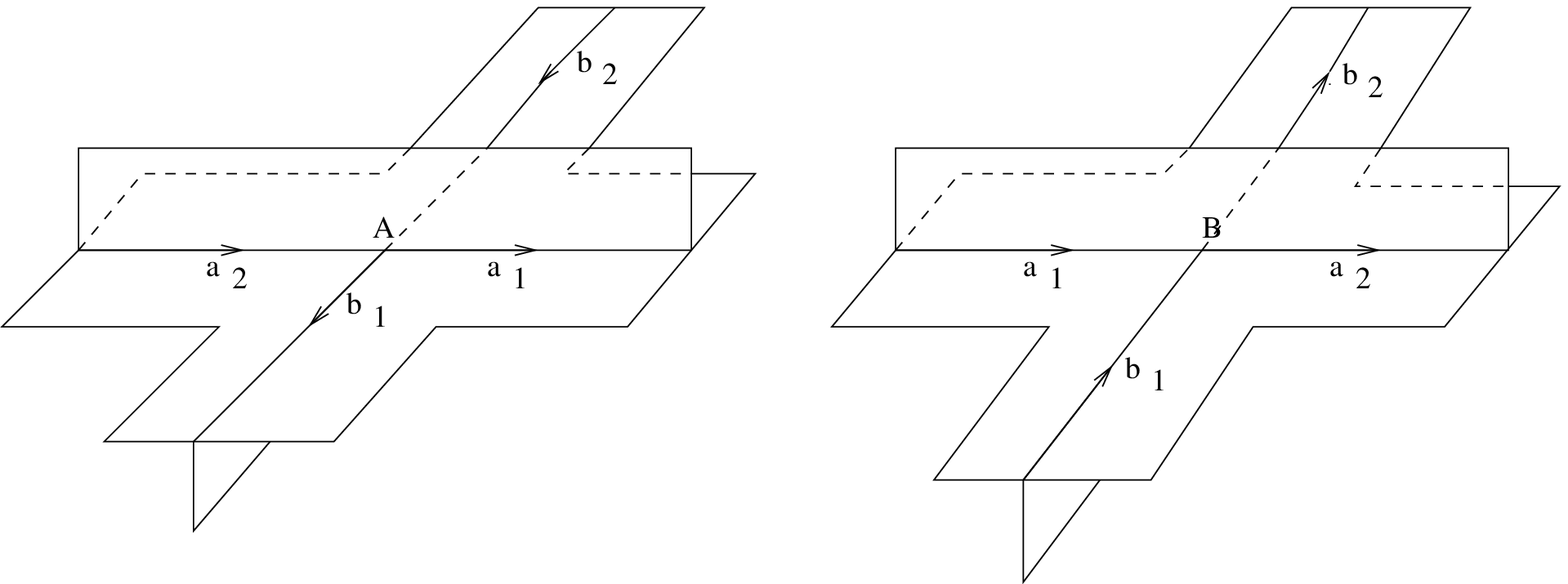}}
\caption{\label{}}
\end{center}
\end{figure}

\begin{tabular}{lllllll}
$(A_{11} A_{12} A_{13})$ & $\sim$ & $(B_{12} B_{11} B_{13})$ & , & $(A_{21} A_{22} A_{23})$ & $\sim$ & $(B_{41} B_{42} B_{43})$  \\[1.ex]

$(A_{31} A_{32} A_{33})$ & $\sim$ & $(B_{32} B_{31} B_{33})$ & , & $(A_{41} A_{42} A_{43})$ & $\sim$ & $(B_{21} B_{22} A_{23})$

\end{tabular} \\[0.3ex]

In this case, the 1-cells $a_{1}$ and $a_{2}$ are even and  the 1-cells $b_{1}$ and $b_{2}$ are odd, $V_{A}\cup V_{B} /\sim$ is embedded into a non-orientable handlebody with an orientable 1-handle containing $b_{1}$, an orientable 1-handle containing $b_{2}$, a non-orientable 1-handle containing $a_{1}$ and a non-orientable 1-handle containing $a_{2}$. Moreover, $Fr(V_{A}\cup V_{B} /\sim)$ consists of four closed paths which are also in the boundary of this handlebody. These four closed paths are homotopic, in  $V_{A}\cup V_{B} /\sim$, to $a_{1}b_{2}a_{1}b_{1}^{-1}$, $a_{1}a_{2}$, $a_{2}b_{1}a_{2}b_{2}^{-1}$ and $b_{1}b_{2}$. We remark that the even 1-cells are present twice in the first three closed paths and absent in the last path. This means, according to our above-mentioned criterion, that the attachement of four disks along the four boundary curves that make up  $Fr(V_{A}\cup V_{B} /\sim)$, produces a standard complex which can be embedded in some orientable 3-manifold.\smallskip

We then attach four disks along these four closed paths and we obtain a polyhedron $\mid {\cal K} \mid$ with the following homotopy relations: \medskip

\begin{tabular}{lllll}

1) $a_{1} b_{2} a_{1} b_{1}^{-1} \sim 1$     & &    2) $a_{1} a_{2} \sim 1$     & & \\[1ex]

3) $a_{2} b_{1} a_{2} b_{2}^{-1} \sim 1$     & &     4) $b_{1} b_{2} \sim 1$     

\end{tabular} \medskip

Taking A as the base point and $a_{1}$ as an auxiliary path leading to vertice B, then the generating path class of the fundamental group of $\mid{\cal K}\mid$ are represented by:\medskip

$A_{1} = a_{1} a_{1}^{-1}$ , $A_{2} = a_{1} a_{2}$ , $B_{1} = b_{1} a_{1}^{-1}$ , $B_{2} = a_{1} b_{2}$ .\smallskip

It follows that $A_{1} \sim 1$ and the homotopy relations can be re-stated by: \medskip

\begin{tabular}{lllll}
1) $B_{1} \sim B_{2}$  &  &   2) $A_{2} \sim 1$  \\[1ex]

 3) $ A_{2} B_{1} A_{2} \sim B_{2}$  &  &  4) $ B_{1} B_{2} \sim 1 $. \\[1ex] 

\end{tabular}

A prompt calculation will then show that $ \pi_{1}(\mid{\cal K}\mid, A )  =  \Z_{2} $.

Let us now consider $\mid{\cal K}\mid$ embedded in an orientable 3-manifold and let $N$ be a regular neighborhood of $\mid{\cal K}\mid$.

The Euler characteristic of $\mid{\cal K}\mid$ is 2, being also the Euler characteristic of $N$ hence this neighborhood is orientable and $\partial N $ is a disjoint union of two spheres of dimension 2. We thus obtain an orientable and closed 3-manifold $M^{3}$ by gluing together two 3-balls on $ \partial N $.

It should be remarked that $\mid{\cal K}\mid$ is embedded in $M^{3}$ but, on the other hand, the  homotopy relations $ a _{1}b_{2}a_{1}b_{1}^{-1} \sim 1$ and $a_{2}b_{1}a_{2}b_{2}^{-1} \sim 1$ define two $M\ddot{o}bius$ strips with the same boundary curve $b_{1} b_{2} $. The $M\ddot{o}bius$ strips are glued together by this boundary curve and form a Klein bottle embedded in $\mid{\cal K}\mid$.\medskip

In fact, $\mid {\cal K} \mid $ is a Klein bottle with two disks $D_{1}$ and $D_{2}$ attached along the paths $a_{1} a_{2}$ and $b_{1} b_{2}$.

We have thus constructed a closed 3-manifold $M^{3}$ whose main features are the following:\medskip

1) Its fundamental group is $\Z_{2} $ and

2) It contains a Klein bottle.

\medskip

Inasmuch, the homotopy relations $ a _{1}b_{2}a_{1}b_{1}^{-1} \sim 1$ and $ b_{1} b_{2} \sim 1 $ define respectively, a $M\ddot{o}bius$ strip and a disk with the same boundary curve $ b_{1} b_{2} $. The $M\ddot{o}bius$ strip and the disk are glued together by their boundary curves and form a projective plane $ RP^{2} $ embedded in $ M^{3} $. A regular neighborhood, $N(RP^{2})$ of $RP^{2}$ in $M^{3}$ is a twisted $[0,1]$-bundle over $RP^{2}$ and $\partial N(RP^{2}) = S^{2} $. The boundary of $M^{3} - Int(N(RP^{2}))$ is $S^{2}$  and by the Van Kampen theorem, it has a trivial fundamental group. Since $M^{3}$ is different from $RP^{3}$, because $RP^{3}$ does not contain a Klein bottle (see [4], [31]), it follows that $M^{3} - Int(N(RP^{2}))$ is not a 3-cell and consequently, $M^{3}$ is reducible.

\medskip

Let $K$ be the Klein bottle, defined by the homotopy relations $ a _{1}b_{2}a_{1}b_{1}^{-1} \sim 1$ and   $a_{2}b_{1}a_{2}b_{2}^{-1} \sim 1$ and let $N(K)$ be a regular neigborhood of $K$. Since $N(K)$ is orientable, $K$ is 1-sided in $N(K)$ and $\partial N(K)$ becomes a torus $T^{2}$. Let $ A_{2} = a_{1} a_{2} $ and $ B_{2} = a_{1} b_{2} $. Then the elements  $A_{2}$ and $B_{2}^{2}$ of the homotopy group $\pi_{1}(K)$ are homotopic in $\pi_{1}(T^{2})$ and  are null homotopic in  $Y^{3} = M^{3} - Int(N(K))$.The neighborhood $N(K)$ is fibered in circles that have homotopy class $B_{2}^{2}$, with two exceptional fibers, each of multiplicity 2, corresponding to the centres of the two $M\ddot{o}bius$ strips of $K$. The orbit surface is a disk. Furthermore, the boundary $\partial Y^{3}$ is equal to $T^{2}$ and $\pi_{1}(Y^{3}) = \Z$. Let us now assume for a moment that $Y^{3}$ is a solid torus. Then, either the fibering extends to $Y^{3}$ and $M^{3}$ becomes a Seifert fibered space or else, $M^{3}$ is the sum of two Lens spaces (see [33]). In the first case, $M^{3}$ must be irreducible (see [8, proposition 1.12] ) and, in the second,  its fundamental group is not $\Z_{2}$. It follows (see [21, Lemma 1]) that $Y^{3}$ is a fake solid torus, that is to say, a fake 3-cell with a 1-handle attached to it along two disjoint disks on its boundary.\medskip

{\bf Proposition 5.5} The space $M^{3}$ is a reducible 3-manifold that satisfies the following properties:

a) There exists, embedded in $M^{3}$, a projective plane $RP^{2}$ such that $M^{3} - Int(N(RP^{2}))$ is a fake 3-cell, where $N(RP^{2})$ is a regular neighborhood of $RP^{2}$ in $M^{3}$.

                                                                          b) There exists, embedded in $M^{3}$, a Klein bottle $K$ such that $M^{3} - Int(N(K))$ is a fake solid torus, where $N(K)$ is a regular neigborhood of $K$ in $M^{3}$. Moreover, $N(K)$ is a Seifert fiber space with two exceptional fibers, each of multiplicity 2, and with orbit surface equal to a disk.\medskip

Finally, we can make use of the universal covering  $W^{3}$ of $M^{3}$ and exhibit a 3-manifold homotopy equivalent to $S^{3}$, though not homeomorphic to it. Let us assume that the only  homotopy 3-sphere is in fact $S^{3}$. Then $W^{3}$  becomes equal to $S^{3}$ and the fundamental group of $M^{3}$, equal to $\Z_{2}$, acts on $S^{3}$. However, the only  element of this group, different from the identity, is an involution of $S^{3}$ and by [13], this involution is conjugate to the antipodal map of $S^{3}$, where after $M^3$ becomes $RP^{3}$. However, $M^{3}$ contains a Klein bottle whereas $RP^{3}$ cannot contain any such bottle. \medskip

{\bf Theorem 5.6} The space $W^{3}$ is a homotopy 3-sphere, non-homeomorphic to $S^{3}$.\medskip

{\bf Remark 1:} A spine of $W^3$ is made up by the universal covering space of $\mid {\cal K} \mid $. This spine is a torus $T^2$ with four disks, two disks attached along two distinct meridians and the other two attached along two distinct parallels. This polyhedron is embeddable in $R^3$. Taking a regular neighborhood of it in $R^3$ and pasting four 3-cells along the four 2-spheres forming the boundary of this regular neighborhood, we obtain $S^3$. This means that two homeomorphic standard complexes can have non-homeomorphic regular neighborhoods, which is not an isolated fact. In general, the standard spine of a homotopy 3-sphere is a polyhedron made up only of orientable surfaces that are glued along their boundary curves or else along closed curves in the interior of the surfaces. This standard spine is embeddable in $R^3$ and is also a spine of $S^3$.  Hence, this does not comply with the theorem 1 and the corollary 1 of [5]. \medskip 

{\bf Remark 2:} Let us take the two vertex neighborhoods $V_{A}$ and $V_{B}$ as are showed in example 5.4 and consider the afterwards indicated matching relations:\medskip

\begin{tabular}{lllllll}
$(A_{11} A_{12} A_{13})$ & $\sim$ & $(B_{32} B_{31} B_{33})$ & , & $(A_{21} A_{22} A_{23})$ & $\sim$ & $(B_{42} B_{41} B_{43})$  \\[1.ex]

$(A_{31} A_{32} A_{33})$ & $\sim$ & $(B_{12} B_{11} B_{13})$ & , & $(A_{41} A_{42} A_{43})$ & $\sim$ & $(B_{22} B_{21} A_{23})$

\end{tabular} \medskip

In this case, the 1-cells $a_{1}$, $a_{2}$, $b_{1}$ and $b_{2}$ are all odd and $V_{A}\cup V_{B} /\sim$ is embedded into an orientable handlebody. Moreover, $Fr(V_{A}\cup V_{B} /\sim)$ consists of four closed paths which are also in the boundary of this handlebody. These four closed paths are homotopic, in  $V_{A}\cup V_{B} /\sim$, to $a_{1}b_{2}a_{2}^{-1}b_{1}^{-1}$, $a_{1}a_{2}$, $a_{1}b_{1}^{-1}a_{2}^{-1}b_{2}$ and $b_{1}b_{2}$. The attachement of four disks along the four boundary curves that make up  $Fr(V_{A}\cup V_{B} /\sim)$, produces a standard complex whose polyhedron is a torus with two discs attached, one along a meridian and the other along a curve that turns twice around the parallels and once around the meridians. It is a spine of $RP^{3}$ and its universal covering is a torus $T^{2}$ with two discs attached along two distinct meridians and two discs attached along two distinct parallels. As we can see, this polyhedron and the standard spine of $W^3$ are the same. This show us that $S^{3}$ and $W^{3}$ have homeomorphic standard spines and the difference between them is very subtle.\medskip

{\bf Remark 3:} We can construct an infinite number of non-homeomorphic homotopy 3-spheres by taking $W^{3}$ and proceeding, by induction, on the connected sums, thus defining a sequence of non-homeomorphic homotopy 3-spheres by setting $M_1$ = $W^{3}$ and $M_{n+1}$ = $M_n$$\#$ $W^{3}$, n = 1, 2, $\ldots$(for the definition of connected sums, see [9]).\medskip

\subsection*{\S6 Thruston Geometrization Conjecture}
There are many conjectures related to that of Poincaré  but we just call the attention upon the  Thurston geometrization conjecture [30]. \smallskip

{\it Every oriented prime closed 3-manifold can be cut along incompressible tori, so that the interior of each of the resulting manifolds has a geometric structure with finite volume. }

\smallskip

We know that if the torus is incompressible, its fundamental group injects into the fundamental group of the 3-manifold. Thus, when the prime closed 3-manifold is a homotopy 3-sphere, the decomposition referred to in the Thruston geometrization conjecture must be  trivial and it implies that the homotopy 3-sphere has just one of the eight geometric structures and the only compact model is $S^{3}$. We thus conclude that this conjecture is not valid as well.\smallskip

{\it Acknowlegments:} I would like to thank Antonio Kumpera for precious help and encouragements.

\end{document}